\documentclass[lettersize,journal]{IEEEtran}
\usepackage{amsmath,amsfonts,amssymb}
\usepackage{algorithmic}
\usepackage{algorithm}
\usepackage{array}
\usepackage[caption=false,font=normalsize,labelfont=sf,textfont=sf]{subfig}
\usepackage{textcomp}
\usepackage{stfloats}
\usepackage{url}
\usepackage{verbatim}
\usepackage{graphicx}
\usepackage{epstopdf}
\usepackage{cite}
\usepackage{mathrsfs}
\usepackage{tabularx}
\usepackage{float}
\usepackage{overpic} 
\usepackage{subfig}
\begin{document}

\title{Fully Distributed Adaptive Nash Equilibrium Seeking Algorithm for Constrained Noncooperative Games with Prescribed Performance}
\author{Sichen Qian
\thanks{Sichen Qian is with the School of Mathematics, Southeast University, Nanjing, 210096, China (email:qiansichen@foxmail.com).}}

\maketitle

\begin{abstract}
This paper investigates a fully distributed adaptive Nash equilibrium (NE) seeking algorithm for constrained noncooperative games with prescribed-time stability. On the one hand, prescribed-time stability for the proposed NE seeking algorithm is obtained by using an adaptive penalty technique, a time-varying control gain and a cosine-related time conversion function, which extends the prior asymptotic stability result. On the other hand, uncoordinated integral adaptive gains are incorporated in order to achieve the fully distribution of the algorithm. Finally, the theoretical result is validated through a numerical simulation based on a standard power market scenario.
\end{abstract}

\section{Introduction}
\IEEEPARstart{N}{oncooperative} games\cite{Nash1950}, an essential aspect of game theory, has been the subject of intense discussion recently due to their unique integration of competitive and collaborative dynamics\cite{Maskery2009, Ershen2024}. 
Also, numerous novel continuous-time algorithms for NE seeking with constraints on players' action profiles has been designed in \cite{Sun2020,Sun2021,Jia2021}. In \cite{Sun2020}, the authors introduce both the exact \( \mathit{l}_1 \) penalty function and the squared \( \mathit{l}_2 \) penalty function with an adjustable penalty parameter which is related to global information and updated based on interaction and search feedback. To overcome the excessive dependence on global information and interactions, the adaptive penalty technique is introduced in \cite{Jia2021} which features a penalty factor adaptively updating according to the constraints violation. To be more specific, the authors in \cite{Jia2021} employ an adaptive penalty gain to transform a constrained distributed resource allocation problem into its unconstrained counterpart which is easier to be dealt with and achieve an asymptotic convergence rate. 

\par
Most crucially, to guarantee the convergence of a closed-loop system, a two-time scale structure is included in \cite{Ye2017,Ye2019,Gadjov2018} so that the consensus part is faster than the optimization part. However, one disadvantage of this technique is that the precise measurement of the singular perturbation parameter is dependent on several centralized information which is difficult to obtain, such as the total number of players and the communication topology between players. But as result of their excessive dependencies on several global information, these approaches lose part of their applicability properties when a new player is added into the game. To deal with the concern, several researchers have put forward some fully distributed NE seeking algorithm (see eg., \cite{Bianchi2020,Ye2021,Ma2022}). Specifically in \cite{Ye2021}, a node-based and an edge-based control law are proposed by adjusting each player's weight on consensus error, and globally asymptotically stability is reached. 
\par
Additionally, the settling time for system states has always been an attractive topic in distributed NE seeking algorithm design. Generally, it can be categorized into finite-time, fixed-time, predefined-time and prescribed-time stability. Compared with the other categories, prescribed-time stability offers a distinct advantage in that the settling time can be predetermined and remains independent of the players' initial states or any global parameters. In recent years, several literature \cite{Tao2022,Feng2020,Zhao2023,Sun2023} has explored the distributed NE seeking with prescribed performance through similar time conversion functions (eg. $T(1-e^{-t})$ in \cite{Feng2020,Tao2022} and $\frac{T}{T-t}$ in \cite{Song2023,Gong2021}). Nevertheless, the existing time conversion technique reveals several disadvantages: On one hand, exponential term decays quickly which might causes rapid changes in system state variables; On the other hand, when approaching the prescribed time $T$, a relatively slow convergence rate may result in poor convergence. Hence, it is desirable to propose a novel cosine-related time conversion function herein which can relief the above concerns.
\par
Motivated by the above studies, this paper presents a fully distributed adaptive prescribed-time NE seeking algorithm for constrained cooperative games. To the best of our knowledge, this specific distributed NE seeking algorithm has not yet been proposed in existing literature. In contrast with previous works in \cite{Sun2020, Sun2021,Jia2021}, a time-varying gain is introduced to further the existing asymptotic results. To achieve fully distribution, uncoordinated integral adaptive gains are introduced while prescribed-time performance is maintained through a cosine-based scheme, which distinguishes from the method utilized in \cite{Feng2020, Tao2022,Song2023,Gong2021}. Generally our approach stands out from existing methods by offering a faster algorithm convergence rate along with fully distribution and adaptive characteristics.
\par
The remainder of this paper is organized as follows. Section \ref{pre} introduces some preliminaries and problem description. Section \ref{main} details the design of our novel NE seeking algorithm and conducts a thorough analysis of prescribed-time convergence. Section \ref{sti} presents a numerical example to testify the theoretical results. Finally, conclusions are drawn in Section \ref{con}.
\par
\textit{Notations:} 
$\mathbb{R}, \mathbb{R}^{m}\text{ and }\mathbb{R}^{m\times n}$ represent the set of real numbers, the set of all $m$ dimension vectors with all component belongs to $\mathbb{R}$ and the set of all $m\times n$ matrices with all component belongs to $\mathbb{R}$ respectively. $[\cdot]^{\top}$ denotes the matrix transpose and $\|\cdot\|$ denotes the Euclidean norm of vector  or the induced 2 norm of matrix respectively.  $\mathbf{0}_{m}=[0,0,\cdots,0]^{\top}\in\mathbb{R}^{m}, \mathbf{1}_{m}=
[1,1,\cdots,1]^{\top}\in\mathbb{R}^{m},\mathbf{I}_{m}$ denotes a $m\times m$ unit matrix. diag$\left\{\xi_{1}, \xi_{2}, \cdots, \xi_{N}\right\}$ denotes a diagonal matrix whose main diagonal element is $\{\xi_{1}, \xi_{2}, \cdots, \xi_{N}\}$. $\otimes$ denotes the Kronecker product. $\lambda_{2}(\cdot), \lambda_{\min}(\cdot)$ represent the smallest positive eigenvalue and the smallest eigenvalue of a matrix respectively. $\partial f(x)$ denotes the subdifferential set of a convex function $f$ at point $x$.
\section{Preliminaries And Problem Description}\label{pre}
\subsection{Graph Theory}
An undirected graph $\mathcal{G:=}\mathcal{(V,E)}$ is made up of a node set $\mathcal{V:=}$\{1,~2,~$\cdots$,~$\mathit{N}$\} and an edge set $\mathcal{E}$. $\mathcal{N}_{i}:=\{\mathit{j:(j,i)}\in$ $\mathcal{E}$\} stands for the neighbors of node $\mathit{i}$, where ($\mathit{j,i})\in\mathcal{E}$ means that node $\mathit{i}$ can communicate with node $\mathit{j}$. A path is a sequence of isolated vertices such that any pair of vertices appearing consecutively is an edge of graph $\mathcal{G}$. A graph $\mathcal{G}$ is said to be connected if there is a path between any two nodes. $\mathcal{A}=\{a_{ij}\}\in\mathbb{R}^{N\times N}$ is the adjacency matrix, where $a_{ii}=0$, $a_{ij}>0$ if $j\in\mathcal{N}_i$ and else, $a_{ij}=0$. Denote $d_i=\sum^N_{j=1}a_{ij}$ as the degree of node $i$. The Laplacian matrix of $\mathcal{G}$ is denoted by $\mathcal{L}=\mathcal{D}-\mathcal{A}$, where $\mathcal{D}=\text{diag}\{d_1,d_2,\cdots,d_N\}$.

\newtheorem{theorem}{\bf Definition}[section]
\subsection{Prescribed-time Stability}
\begin{theorem}\label{thm3}\cite{Zhu2023}
Consider a continuous time system
\begin{equation}\label{xitong}
	\dot{x}(t)=T(t, T_{p})f(x(t))
\end{equation}
where time $t>0, T(t, T_{p}): \mathbb{R}^{+}\cup{\{0\}}\times\mathbb{R}^{+}\rightarrow\mathbb{R}$ is a function with time $t$ as its independent variable, and $T_{p}$ is the prescribed time.  $f:\mathbb{R}^{n}\rightarrow\mathbb{R}$ is a mapping with $x$ as its independent variable. 
Say that the system \eqref{xitong} achieves prescribed-time stable to the point $x^{*}$ at $T _ {p}$, if for any $x(0)$, $\lim\nolimits_{t\rightarrow T_{p}^{-}}\vert\vert x(t)-x^{*}\vert\vert=0$ holds. 
\end{theorem}
\par
In an attempt to achieve the prescribed performance, a time conversion function $k(t, T_{p}):[0, T_{p})\longmapsto \mathbb{R}^{+}$ related to cosine function is employed as
\begin{equation}\label{trans}
	k(t)=\sec^{2}(\frac{\pi}{2}\frac{t}{T_{p}}), t\in[0, T_{p})
\end{equation}
where $T_{p}>0$ is the prescribed time appointed beforehand. And the conversion function $k(t)$ has the properties that $k(0)=1$ and $\lim_{t\rightarrow T_{p}^{-}}k(t)=+\infty$.

\subsection{Problem Description}
Consider a noncooperative game $\Lambda\triangleq\{\mathcal{V}, J_{i}, \Omega_{i}\}$ with $N$ players on an undirected connected graph $\mathcal{G}=(\mathcal{V},\mathcal{E})$, where $\mathcal{V}=\{1, 2, \cdots, N\}.$ And each player dominates his own action profile $x_{i}\in\Omega_{i}$ where the convex set $\Omega_{i}=\{x_{i}\in\mathbb{R}^{n_{i}}\vert g_{i}(x_{i})\leqslant 0\}$ is the feasible action profile set for the $i$th player. Define $x_{-i}=[
(x_{1})^{\top},\cdots,(x_{i-1})^{\top},(x_{i+1})^{\top},\cdots,(x_{N})^{\top}
]^{\top}$. And the $i$th player's cost function in the noncooperative game $\Lambda$ is denoted as $J_{i}(x_{i}, x_{-i}):\mathbb{R}^{n}\rightarrow\mathbb{R}, n=\sum_{i\in\mathcal{V}}n_{i}$ which is a continuously differentiable convex function with respect to the variable $x_{i}$. Denote the vector whose components are all players' own strategies as $x=[
x_{1},\cdots, x_{i},\cdots,x_{N}
]^{\top}$. 
\par 
Specifically, the $i$th player's goal is shown as follows:
	\begin{equation*}
		\begin{aligned}
			\text{minimize} \quad &J_{i}(x_{i}, x_{-i}) \\
			\text{subject to} \quad&x_{i}\in\Omega_{i}.
		\end{aligned}
	\end{equation*}
\par 
For noncooperative game $\Lambda$, NE is a preferred optimal solution defined in Definition \ref{ne}.
\begin{theorem}\label{ne}\cite{Nash1950}
An action profile $x_{i}^{*}=(x_{i}^{*}, x_{-i}^{*})\in\Omega, \Omega=\Omega_{1}\times\Omega_{2}\times\cdots\times\Omega_{N}\subseteq\mathbb{R}^{Nn}$ is said to be a Nash equilibrium (NE) point for the noncooperative game $\Lambda$ under the condition that
\begin{equation*}
	J_{i}(x_{i}^{*}, x_{-i}^{*})\leqslant J_{i}(x_{i}, x_{-i}^{*}), \forall x_{i}\in\Omega_{i}, \forall i \in \mathcal{V}
\end{equation*}
where $x_{-i}^{*}=[x_{1}^{*}, \cdots, x_{i-1}^{*}, x_{i+1}^{*}, \cdots, x_{N}^{*}]^{\top}, x_{i}^{*}$ is the $i$th player's Nash equilibrium (NE) action profile.
\end{theorem}

\par
To guarantee the existence and uniqueness of NE for noncooperative game $\Lambda$, the following assumption is frequently referred to.
\newtheorem{lem}{\bf Lemma}[section]

\newtheorem{assum}{\bf Assumption}[section]
\begin{assum}\label{gm}
The pseudo-gradient for $J_{i}(x_{i}, x_{-i})$ is $\mu$-strongly monotone and $l$-Lipschitz continuous.
\end{assum}
\section{Main Results}\label{main}
\subsection{Fully Distributed Adaptive NE Seeking}
A fully distributed adaptive NE seeking algorithm with a special design of prescribed-time stability is proposed with passivity increment as in \eqref{algo_1}. Firstly, in distributed computing, the leader-follower estimation approach is necessary since each player's cost function requires knowledge on all other players' actions, in addition to their own. Thus we define the $i$th player's estimation on the action profiles of all other players as $x^{i}=[
(x^{i}_{1})^{\top},\cdots,(x^{i}_{N})^{\top}
]^{\top}$ where $x^{i}_{j}$ is the $i$th player's estimation on the $j$th player's strategy. 

\begin{equation}\label{algo_1}	\left\{
	\begin{aligned}
		\dot{x}_{i}	=&-k(t)[\nabla_{i}J_{i}(x_{i}, x^{i}_{-i})+\sigma_{i}\eta_{i}+\mathcal{R}_{i}\sum\limits_{j\in\mathcal{N}_{i}}(\omega_{i}\rho^{i}-\omega_{j}\rho^{j})],\\
		\dot{x}^{i}_{-i}=&-k(t)\mathcal{S}_{i}\sum_{j\in\mathcal{N}_{i}}(\omega_{i}\rho^{i}-\omega_{j}\rho^{j}),\\
		\dot{\sigma}_{i}=&k(t)q(t)G_{i}(x_{i}),\\
			\dot{\omega}_{i}=&k(t)q(t)\gamma_{i}\Vert \rho^{i}\Vert^{2},\\
			\rho^{i}=&\sum_{j\in\mathcal{N}_{i}}(x^{i}-x^{j})
	\end{aligned}
	\right.
\end{equation}
where $k(t)$ is as defined in \eqref{trans}. And $\dot{q}(t)=k(t)q(t)$ with $q(0)>0$,  $\eta_{i}\in\partial_{i}G_{i}(x_{i}(t))$. Additionally, the matrix $\mathcal{R}_{i}$ and the matrix $\mathcal{S}_{i}$ are respectively defined as
\begin{equation*}
	\mathcal{R}_{i}=
	\begin{bmatrix}
		\mathbf{0}_{n_{i}\times n_{<i}} & \mathbf{I}_{n_{i}} &\mathbf{0}_{n_{i}\times n_{>i}}
	\end{bmatrix},
\end{equation*}
\begin{equation*}
	\mathcal{S}_{i}=
	\begin{bmatrix}
		\mathbf{I}_{n_{<i}}&\mathbf{0}_{n_{<i}\times n_{i}}&\mathbf{0}_{n_{<i}\times n_{>i}}\\
		\mathbf{0}_{n_{>i}\times n_{<i}}&\mathbf{0}_{n_{>i}\times n_{i}}&\mathbf{I}_{n_{>i}}
	\end{bmatrix}
\end{equation*}
where $n_{<i}=\sum_{j<i, i, j\in\mathcal{V}}n_{j}, n_{>i}=\sum_{j>i, i, j\in\mathcal{V}}n_{j}.$
 And the designed algorithm \eqref{algo_1} can be further rewritten into the compact form as is shown in \eqref{algo_2}.
\begin{equation}\label{algo_2}\left\{
	\begin{aligned}
		\dot{X}&=-k(t)\left\{\mathcal{R}^{\top}\left[\mathbf{F}(X)+\eta(\mathbf{x})\sigma\right]+\mathbf{L}Z\rho\right\}\\
		\dot{\sigma}&=k(t)q(t)G(\mathbf{x})\\
		\dot{\omega}&=k(t)q(t)D(\rho)^{\top}(\Gamma\otimes\mathbf{I}_{n})\rho\\
		\rho&=\mathbf{L}X
	\end{aligned}\right. 
\end{equation}
where
\begin{flalign*}
	&& X & = \left[(x^{1})^{\top}, (x^{2})^{\top}, \cdots, (x^{N})^{\top}\right]^{\top} & \\
	&& \mathbf{F}(x) & = \left[\nabla_{1}J_{1}(x^{1}), \nabla_{2}J_{2}(x^{2}), \cdots, \nabla_{N}J_{N}(x^{N})\right]^{\top}  & \\
	&& \mathcal{R} & = \text{diag}\{\mathcal{R}_{1}, \mathcal{R}_{2}, \cdots, \mathcal{R}_{N}\} & \\
	&& \eta& = \text{diag}\{\eta_{1}, \eta_{2}, \cdots, \eta_{N}\} &\\
	&&\sigma &={[\sigma_{1}, \sigma_{2}, \cdots, \sigma_{N}]}^{\top}&\\
	&&\mathbf{L} &=\mathcal{L} \otimes \mathbf{I}_{n}&\\
	&&G(\mathbf{x})&=[G_{1}(x_{1}), G_{2}(x_{2}), \cdots, G_{N}(x_{N})]^{\top}&\\
	&&\omega&=[\omega_{1}, \omega_{2},\cdots,\omega_{N}]&\\
	&&\rho&=[\rho^{1},\rho^{2},\cdots, \rho^{N}]&\\
	&&Z&=\text{diag}\{\omega_{1}\mathbf{I}_{n},\omega_{2}\mathbf{I}_{n},\cdots,\omega_{N}\mathbf{I}_{n}\}&\\
	&&\Gamma&=\text{diag}\{\gamma_{1},\gamma_{2}, \cdots, \gamma_{N}\}&\\
	&&D(\rho)&=\text{diag}\{\rho^{1},\rho^{2},\cdots,\rho^{N}\}\subseteq\mathbb{R}^{Nn\times N}&.
\end{flalign*}
\par
We first prove that  $\widetilde{X}$ satisfying $\dot{\widetilde{x}}=\mathbf{0}_{Nn}, \dot{\widetilde{\sigma}}=\mathbf{0}_{N},\dot{\widetilde{\omega}}=\mathbf{0}_{N}$ is the NE point $\mathbf{1}_{N}\otimes x^{*}$ for noncooperative game $\Lambda$.
\newtheorem{dingli}{\bf Theorem}[section]
\begin{dingli}\label{dengjia}
 $\widetilde{X}$ satisfying $\dot{\widetilde{X}}=\mathbf{0}_{Nn},\dot{\widetilde{\sigma}}=\mathbf{0}_{N},\dot{\widetilde{\omega}}=\mathbf{0}_{N}$, is equivalent to the NE point $x^{*}$ of noncooperative game $\Lambda$ and $\widetilde{X}^{i}=\widetilde{X}^{j}=x^{*}, \forall i, j \in \mathcal{V}$ holds.
\end{dingli}
\begin{IEEEproof}
Noticing that $(\widetilde{X}, \widetilde{\sigma},\widetilde{\omega})$ satisfies $\dot{\widetilde{X}}=\mathbf{0}_{Nn},\dot{\widetilde{\sigma}}=\mathbf{0}_{N},\dot{\widetilde{\omega}}=\mathbf{0}_{N}$ respectively and for $t\in[0,T_{p}), k(t)>0, q(t)>0$, then it holds that 
\begin{align}
	\mathbf{0}_{Nn}&=-\mathcal{R}^{\top}(\mathbf{F}(\widetilde{X})+\eta(\widetilde{\mathbf{x}})\widetilde{\sigma})-\mathbf{L}\widetilde{Z}\widetilde{\rho},\label{equil}\\
	\mathbf{0}_{N}&=G(\widetilde{\mathbf{x}}),\\
	\mathbf{0}_{N}&=D(\widetilde{\rho})^{\top}(\Gamma\otimes\mathbf{I}_{n})\widetilde{\rho},\label{er}\\
	\widetilde{\rho}&=\mathbf{L}\widetilde{X},\nonumber\\
	\widetilde{Z}&=\text{diag}\left\{\widetilde{\omega}_{1}\mathbf{I}_{n}, \widetilde{\omega}_{2}\mathbf{I}_{n},\cdots,\widetilde{\omega}_{N}\mathbf{I}_{n}\right\}\nonumber
\end{align}
where $\widetilde{\mathbf{x}}=\left[(\widetilde{x}_{1})^{\top}, \cdots, (\widetilde{x}_{i})^{\top}, \cdots, (\widetilde{x}_{N})^{\top}\right]^{\top}$.
Noticing \eqref{er} holds and $\gamma_{i}>0, i\in\left\{1,2,\cdots.N\right\}$, it can be deduced that $\mathbf{L}\widetilde{X}=\mathbf{0}_{Nn}$, which further indicates that $\widetilde{X}^{i}=\widetilde{X}^{j}=x^{*}, \forall i, j \in \mathcal{V}$. 
\par
On the other hand, let $\widetilde{X}=\mathbf{1}_{N}\otimes x^{*}, x^{*}\in\Omega$. 
By multiplying $\mathbf{1}_{N}^{\top}\otimes\mathbf{I}_{n}$ to the left at both ends of \eqref{equil} and based on the identity equation $\mathcal{L}\mathbf{1}_{N}=\mathbf{0}_{N}$, it can be deduced that 
$
	\mathbf{0}_{Nn}=-\mathcal{R}^{\top}(\mathbf{F}(\widetilde{X})+\eta(\widetilde{\mathbf{x}})\widetilde{\sigma}).
$
By noticing that $\mathbf{L}\widetilde{X}=\mathbf{0}_{Nn}$, it can be obtained that $\dot{\widetilde{X}}=\mathbf{0}_{Nn},\dot{\widetilde{\sigma}}=\mathbf{0}_{N}$ and $\dot{\widetilde{\omega}}=\mathbf{0}_{N}$. In conclusion, Lemma \ref{dengjia} holds.
\end{IEEEproof}
\par
Next,  two arguments are given to facilitate the proof for Theorem \ref{the}
\begin{lem}\label{xzc1}\cite{Cui2021}
		Consider the following system
		\begin{equation}\label{Vxt}
			\dot{V}(t)=-\theta k(t)V(t), V(0)>0, \theta>0
		\end{equation} 
		where $k(t)$ is as defined in \eqref{trans}, then the system \eqref{Vxt} is prescribed-time stable to the origin at time $T_{p}$.
\end{lem}
\par
The detailed proof for Lemma \ref{xzc1} are similar to the proof for Lemma 1 in \cite{Cui2021}, and thus it is omitted therein.

\newtheorem{col}{\bf Corollary}
\begin{col}\label{xzc2}
Suppose there exists a Lyapunov function $V(x)$ satisfying the condition that $\exists M>0, V(x)\geqslant M\Vert x\Vert^{2}$ and \begin{equation}
	\dot{V}(x(t))\leqslant-\theta k(t) V(x(t)),V(x(0))>0, \theta>0
\end{equation}
where $k(t)$ is as defined in \eqref{trans}, then the according system
\begin{equation}\label{ss}
	\dot{x}(t)=k(t)f(x(t))
\end{equation}
is prescribed-time stable to the origin at time $T_{p}$.
\end{col}
\begin{IEEEproof}
	From the proof for Lemma \ref{xzc1}, the inequality holds that
	\begin{align}\label{transs}
	0\leqslant\vert\vert x(t) \vert\vert^{2}&\leqslant M^{-1}\exp^{-\theta\int_{0}^{t}k(\tau)d\tau}V(0)\nonumber\\
	&\leqslant M^{-1}\exp^{-\frac{2\theta T_{p}}{\pi}(\alpha(t)-\alpha(0))}V(0)
	\end{align} where $\alpha(t)=\tan(\frac{\pi}{2}\frac{t}{T_{p}})$ and $\alpha(0)=0, \lim_{t\rightarrow T_{p}^{-}}\alpha(t)=+\infty$ hold. Then by noticing that initial value $x(0)$ and $V(0)$ are all positive, it can be deduced that $\lim_{t\rightarrow T_{p}^{-}}x(t)=0$ according to limit pinch principal.
\end{IEEEproof}

\par
Next the main theorem of this paper is given.
\begin{dingli}\label{the}
	Suppose Assumption \ref{gm} holds and all player's action profile is updated according to \eqref{algo_2}, then all player's action profile is prescribed-time stable to $\mathbf{1}_{N}\otimes x^{*}$ at time $T_{p}$, that is to say, action profile components of all players is prescribed-time stable to the NE point $x^{*}$ of the noncooperative game $\Lambda$ at time $T_{p}$.
\end{dingli}
\begin{IEEEproof}
Firstly, a method for orthogonal decomposition of space $\mathbb{R}^{Nn}$ is introduced.
$\mathbb{R}^{Nn}=\mathbb{C}_{N}^{n}\oplus \mathbb{E}_{N}^{n}$ where $\mathbb{C}_{N}^{n}=\left\{\mathbf{1}_{N}\otimes x\vert x\in\mathbb{R}^{n}\right\}$ is the subspace that satisfies the consensus protocol and $\mathbb{E}_{N}^{n}$ is its corresponding orthogonal complement space. And the following two projection matrices can be defined:
\begin{equation*}
	P_{\mathbb{C}}=\frac{1}{N}\mathbf{1}_{N}\otimes\mathbf{1}^{\top}_{N}\otimes \mathbf{I}_{n}\text{ and }P_{\mathbb{E}}=\mathbf{I}_{Nn}-\frac{1}{N}\mathbf{1}_{N}\otimes\mathbf{1}_{N}^{\top}\otimes\mathbf{I}_{n}.
\end{equation*}
Thus, $X\in\mathbb{R}^{Nn}$ can be decomposed into 
$
X=X^{\parallel}+X^{\perp}
$
where $X^{\parallel}=P_{\mathbb{C}}X\in\mathbb{C}_{N}^{n}$, $ X^{\perp}=P_{\mathbb{E}}X\in\mathbb{E}_{N}^{n}$. 
\par
Secondly, consider a Lyapunov function candidate
\begin{equation*}
\begin{aligned}
V=&\frac{1}{2}\Vert X-\widetilde{X}\Vert^{2}+\frac{1}{2q(t)}\Vert\sigma-\widetilde{\sigma}\Vert^{2}+\frac{1}{2q(t)}\Vert \omega-\widetilde{\omega}\Vert^{2}_{\Gamma^{-1}}\\
=&\frac{1}{2}\Vert X-\widetilde{X}\Vert^{2}+\frac{1}{2q(t)}(\Vert\sigma-\widetilde{\sigma}\Vert^{2}+(\omega-\widetilde{\omega})^{\top}\Gamma^{-1}(\omega-\widetilde{\omega}))
\end{aligned}
\end{equation*}
where $\widetilde{X}=\mathbf{1}_{N}\otimes x^{*}$ and $\widetilde{\omega}$ such that $\omega^{*}:=\min(\widetilde{\omega})>\underline{\omega}= \frac{l^{2}+l\mu}{\mu\lambda_{2}(\mathcal{L})^{2}}$. And since $\dot{q}(t)=k(t)q(t), q(0)>0$, it holds that $q(t)>0, \forall t>0$. Consequently, the above Lyapunov candidate function satisfies that 
$
	V\geqslant\frac{1}{2}\left\|X-\widetilde{X}\right\|^{2}.
$
\par
Next, taking the derivative of $V$ along time $t$ and based on \eqref{algo_2}, it can be deduced that \begin{align}
	\dot{V}=&-k(t)\left(X-\widetilde{X}\right)^{\top}\mathcal{R}^{\top}\left[\mathbf{F}(X)-\mathbf{F}(\widetilde{X})\right]\nonumber\\
	&-k(t)\left(X-\widetilde{X}\right)^{\top}\mathcal{R}^{\top}\left(\eta(\mathbf{x})\sigma-\eta(\widetilde{\mathbf{x}})\widetilde{\sigma}\right)\nonumber\\
	&+\frac{1}{q(t)}k(t)\left(\sigma-\widetilde{\sigma}\right)^{\top}q(t)G(\mathbf{x})\nonumber\\
	&-\frac{\dot{q}(t)}{2q^{2}(t)}\left\|\sigma-\widetilde{\sigma}\right\|^{2}\nonumber\\
	&-k(t)\left(X-\widetilde{X}\right)^{\top}\mathbf{L}Z\mathbf{L}(X-\widetilde{X})\nonumber\\
	&+\frac{1}{q(t)}k(t)(\omega-\widetilde{\omega})^{\top}\Gamma^{-1}D(\rho)^{\top}(\Gamma\otimes\mathbf{I}_{N})\rho\nonumber\\
	&-\frac{\dot{q}(t)}{2q^{2}(t)}\Vert\omega-\widetilde{\omega}\Vert^{2}_{\Gamma^{-1}}\label{Vdaoshu_2}
\end{align}
For the sake of convenience, $\dot{V}$ is separated into $\dot{V}_{1}, \dot{V}_{2}$ and $\dot{V}_{3}$, where $\dot{V}_{1}$ denotes the first addend in \eqref{Vdaoshu_2} , $\dot{V}_{3}$ denotes the last three addends in \eqref{Vdaoshu_2} and $\dot{V}_{2}$ represents the remaining addends.
\par
As for $\dot{V}_{1}$, noting that $X=X^{\parallel}+X^{\perp},\widetilde{X}=\mathbf{1}_{N}\otimes x^{*}$ and $\mathbf{L}X^{\parallel}=0$, it can be deduced that 
\begin{align}
	\dot{V}_{1}=&-k(t)\left(X^{\perp}\right)^{\top}\mathcal{R}^{\top}\left[\mathbf{F}(X)-\mathbf{F}(X^{\parallel})\right]\nonumber\\
	&-k(t)\left(X^{\perp}\right)^{\top}\mathcal{R}^{\top}\left[\mathbf{F}(X^{\parallel})-\mathbf{F}(\widetilde{X})\right]\nonumber\\
	&-k(t)\left(X^{\parallel}-\mathbf{1}_{N}\otimes x^{*}\right)^{\top}\mathcal{R}^{\top}\left[\mathbf{F}(X)-\mathbf{F}(X^{\parallel})\right]\nonumber\\
	&-k(t)\left(X^{\parallel}-\mathbf{1}_{N}\otimes x^{*}\right)^{\top}\mathcal{R}^{\top}\left[\mathbf{F}(X^{\parallel})-\mathbf{F}(\widetilde{X})\right]\nonumber.
\end{align}
With properties of a real symmetric positive definite matrices,   it holds that 
\begin{align}
	\dot{V}_{1}&\leqslant k(t)\Vert X^{\perp}\Vert\Vert\mathcal{R}\Vert\Vert\mathbf{F}(X)-\mathbf{F}(X^{\parallel})\Vert\nonumber\\
	&-k(t)\left(X^{\perp}\right)^{\top}\mathcal{R}^{\top}\left[\mathbf{F}(\mathbf{1}_{N}\otimes x)-\mathbf{F}(\mathbf{1}_{N}\otimes x^{*})\right]\nonumber\\
	&-k(t)\left(X^{\parallel}-\mathbf{1}_{N}\otimes x^{*}\right)^{\top}\mathcal{R}^{\top}\left[\mathbf{F}(X)-\mathbf{F}(X^{\parallel})\right]\nonumber\\
	&-k(t)\left(X^{\parallel}-\mathbf{1}_{N}\otimes x^{*}\right)^{\top}\mathcal{R}^{\top}\left[\mathbf{F}(\mathbf{1}_{N}\otimes x)-\mathbf{F}(\mathbf{1}_{N}\otimes x^{*})\right]\nonumber.
\end{align}
By further referring to Assumption \ref{gm} and the definition of matrix $\mathcal{R}$, it holds obviously that \begin{align}
	\dot{V}_{1}\leqslant&k(t)l\Vert X^{\perp}\Vert^{2}\nonumber\\
	&-k(t)\left(X^{\perp}\right)^{\top}\mathcal{R}^{\top}\left[\mathbf{F}(\mathbf{1}_{N}\otimes x)-\mathbf{F}(\mathbf{1}_{N}\otimes x^{*})\right]\nonumber\\
	&-k(t)\left(x-x^{*}\right)^{\top}\left[\mathbf{F}(X)-\mathbf{F}(X^{\parallel})\right]\nonumber\\
	&-k(t)\left(x-x^{*}\right)^{\top}\left[\mathbf{F}(\mathbf{1}_{N}\otimes x)-\mathbf{F}(\mathbf{1}_{N}\otimes x^{*})\right].\nonumber
\end{align}
Similarly, by employing Assumption \ref{gm} and the definition of matrix $\mathcal{R}$ along with $\Vert x-x^{*}\Vert=\frac{1}{\sqrt{N}}\Vert X^{\parallel}-\widetilde{X}\Vert$, it can be deduced that \begin{align}
	\dot{V}_{1}\leqslant&-k(t)\left[-l\Vert X^{\perp}\Vert^{2}-\frac{2l}{\sqrt{N}}\Vert X^{\perp}\Vert\Vert X^{\parallel}-\widetilde{X}\Vert\right.\nonumber\\
	&+\left.\frac{\mu}{N}\Vert X^{\parallel}-\widetilde{X}\Vert^{2}\right].\label{V1}
\end{align}
As for $\dot{V}_{2}$, 
with properties of a convex function, it holds that 
$
	G(\mathbf{x})\geqslant G(\widetilde{\mathbf{x}})+\left\langle\mathbf{x}-\widetilde{\mathbf{x}}, \eta(\widetilde{\mathbf{x}})\right\rangle,
	G(\widetilde{\mathbf{x}})\geqslant G(\mathbf{x})+\left\langle\widetilde{\mathbf{x}}-\mathbf{x}, \eta(\mathbf{x})\right\rangle.
$
Then it can be deduced that 
\begin{align}
	\dot{V}_{2}\leqslant&k(t)\left[G(\widetilde{\mathbf{x}})-G(\mathbf{x})\right]\sigma+k(t)\left[G(\mathbf{x})-G(\widetilde{\mathbf{x}})\right]\widetilde{\sigma}\nonumber\\
	&+k(t)\left(\sigma-\widetilde{\sigma}\right)^{\top}G(\mathbf{x})-k(t)\frac{1}{2q(t)}\left\|\sigma-\widetilde{\sigma}\right\|^{2}\nonumber\\
	=&k(t)G(\widetilde{\mathbf{x}})\left(\sigma-\widetilde{\sigma}\right)-k(t)\frac{1}{2q(t)}\Vert\sigma-\widetilde{\sigma}\Vert^{2}\nonumber\\
	=&-k(t)\frac{1}{2q(t)}\Vert\sigma-\widetilde{\sigma}\Vert^{2}.\label{VV2_}
\end{align}
As for $\dot{V}_{3}$, we first address the last second addend in \eqref{Vdaoshu_2} as 
\begin{equation}\label{V3}
	\begin{aligned}
&k(t)(\omega-\widetilde{\omega})^{\top}\Gamma^{-1}D(\rho)^{\top}(\Gamma\otimes\mathbf{I}_{N})\rho\\
=&k(t)\sum_{i=1}^{N}(\omega_{i}-\widetilde{\omega}_{i})\rho^{i^{\top}}\rho^{i}\\
=&k(t)\rho^{\top}(Z-\bar{Z})\rho\\
=&k(t)X^{\top}\mathbf{L}(Z-\bar{Z})\mathbf{L}X\\
=&k(t)(X-\widetilde{X})^{\top}\mathbf{L}(Z-\bar{Z})\mathbf{L}(X-\widetilde{X}).
\end{aligned}
\end{equation}
Then, $\dot{V}_{3}$ can be formulated as 
\begin{equation}
	\begin{aligned}
		\dot{V}_{3}=&-k(t)(X-\widetilde{X})^{\top}\mathbf{L}\bar{Z}\mathbf{L}(X-\widetilde{X})-\frac{k(t)}{2q(t)}\Vert\omega-\widetilde{\omega}\Vert^{2}_{\Gamma^{-1}}\\
		\leqslant&-k(t)(X-\widetilde{X})^{\top}\mathbf{L}Z^{*}\mathbf{L}(X-\widetilde{X})-\frac{k(t)}{2q(t)}\Vert\omega-\widetilde{\omega}\Vert^{2}_{\Gamma^{-1}}\\
		=&-k(t)(X^{\parallel}+X^{\perp}-\widetilde{X})^{\top}\mathbf{L}Z^{*}\mathbf{L}(X^{\parallel}+X^{\perp}-\widetilde{X})\\
		&-\frac{k(t)}{2q(t)}\Vert\omega-\widetilde{\omega}\Vert^{2}_{\Gamma^{-1}}\\
		=&-k(t)X^{\perp^{\top}}\mathbf{L}Z^{*}\mathbf{L}X^{\perp}-\frac{k(t)}{2q(t)}\Vert\omega-\widetilde{\omega}\Vert^{2}_{\Gamma^{-1}}\\
		\leqslant&-k(t)\omega^{*}\lambda_{2}(\mathcal{L})^{2}\Vert X^{\perp}\Vert^{2}-\frac{k(t)}{2q(t)}\Vert\omega-\widetilde{\omega}\Vert^{2}_{\Gamma^{-1}}
	\end{aligned}
\end{equation}
where $Z^{*}=\omega^{*}\mathbf{I}_{Nn}$. 
\par
Thus it holds that 
$
	\dot{V}=\dot{V}_{1}+\dot{V}_{2}+\dot{V}_{3}\leqslant-k(t)\Theta\Xi\Theta^{\top}
$
where 
\begin{equation*}
	\Xi=\begin{bmatrix}
		\frac{\mu}{N}&-\frac{l}{\sqrt{N}}&0&0\\
		-\frac{l}{\sqrt{N}}&\omega^{*}\lambda_{2}(\mathcal{L})^{2}-l&0&0\\
		0&0&1&0\\
		0&0&0&1\\
	\end{bmatrix}
\end{equation*}
and $	\Theta=$
\begin{equation*}
	[
\Vert X^{\parallel}-\widetilde{X}\Vert,\Vert X^{\perp}\Vert,\frac{1}{\sqrt{2q(t)}}\Vert\sigma-\widetilde{\sigma}\Vert,\frac{1}{\sqrt{2q(t)}}\Vert\omega-\widetilde{\omega}\Vert_{\sqrt{\Gamma^{-1}}}].
\end{equation*}
Noticing that the matrix $\Xi$ is a symmetric positive definite matrix, it can be reached that 
$
	\dot{V}\leqslant-k(t)\lambda_{\min}(\Xi)V.
$
By letting $M=0.5, \theta=\lambda_{\min}(\Xi)$ in Corollary \ref{xzc2}, it is obvious that the system \eqref{algo_2} is prescribed-time stable.
\end{IEEEproof}
\section{Numerical Stimulation}\label{sti}
\begin{figure}[htbp]
\centering
\includegraphics[scale=0.15]{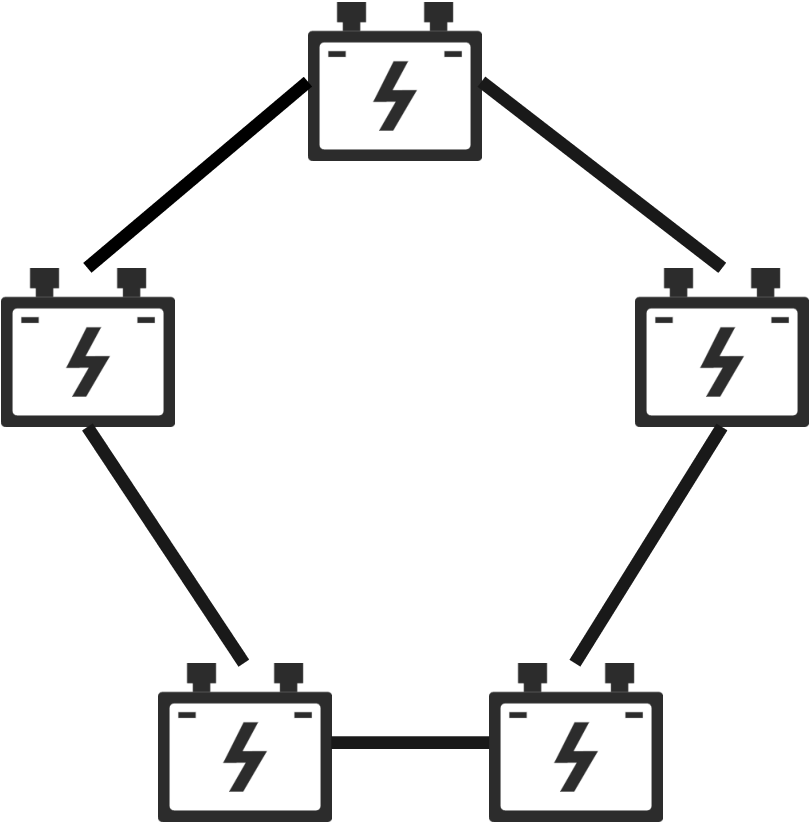}
\caption{Communication Topology Among Five Power Plants}\label{top}
\end{figure}

In this section, a simplified electricity demand response management model is introduced to verify the theoretical result. Consider an electricity consumption game consisting of five power plants, whose electrical appliances can be divided into low, medium and high grades. Specifically the electricity consumption of the $i$th power plant can be expressed as a vector $x_{i}=[x_{i_{1}},x_{i_{2}},x_{i_{3}}]^{\top}$. The communication topology between power plants is depicted in Fig. \ref{top}. 
The energy cost function is set as $J_{i}=a_{i}x^{2}_{i_{1}}+b_{i}x^{2}_{i_{2}}+c_{i}x^{2}_{i_{3}}+\sum_{p,q\in\mathcal{V}\setminus\{i\}}x^{i}_{p}x^{i}_{q}-d_{i}x_{i_{1}}-e_{i}x_{i_{2}}-f_{i}x_{i_{3}}$, and the  inequality constraint function is set as $g_{i}=\alpha_{i}(x_{i_{1}}+x_{i_{2}}+x_{i_{3}})-\beta_{i}$. 

Also several parameters are set as  $T_{p}=10$s, $q(0)=0.001$, $\sigma_{1}(0)=\cdots=\sigma_{5}(0)=50$, $\gamma_{1}=\cdots=\gamma_{5}=1$ and $\omega_{1}(0)=\cdots=\omega_{5}(0)=0.001$. The initial value for each power plant's action profile is preset randomly. 

\begin{figure}[htbp]
\centering
\includegraphics[width=\linewidth]{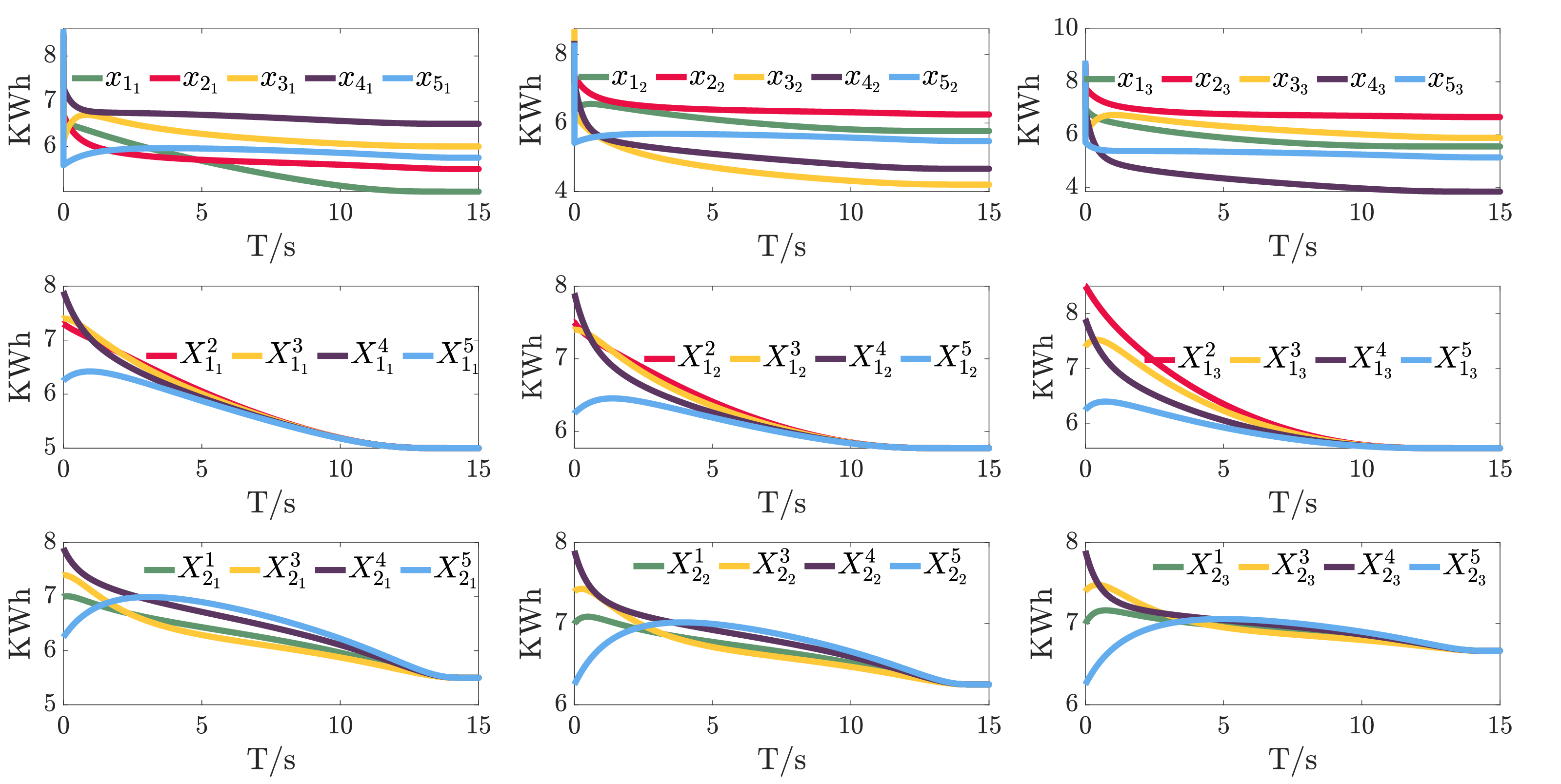}
\includegraphics[width=\linewidth]{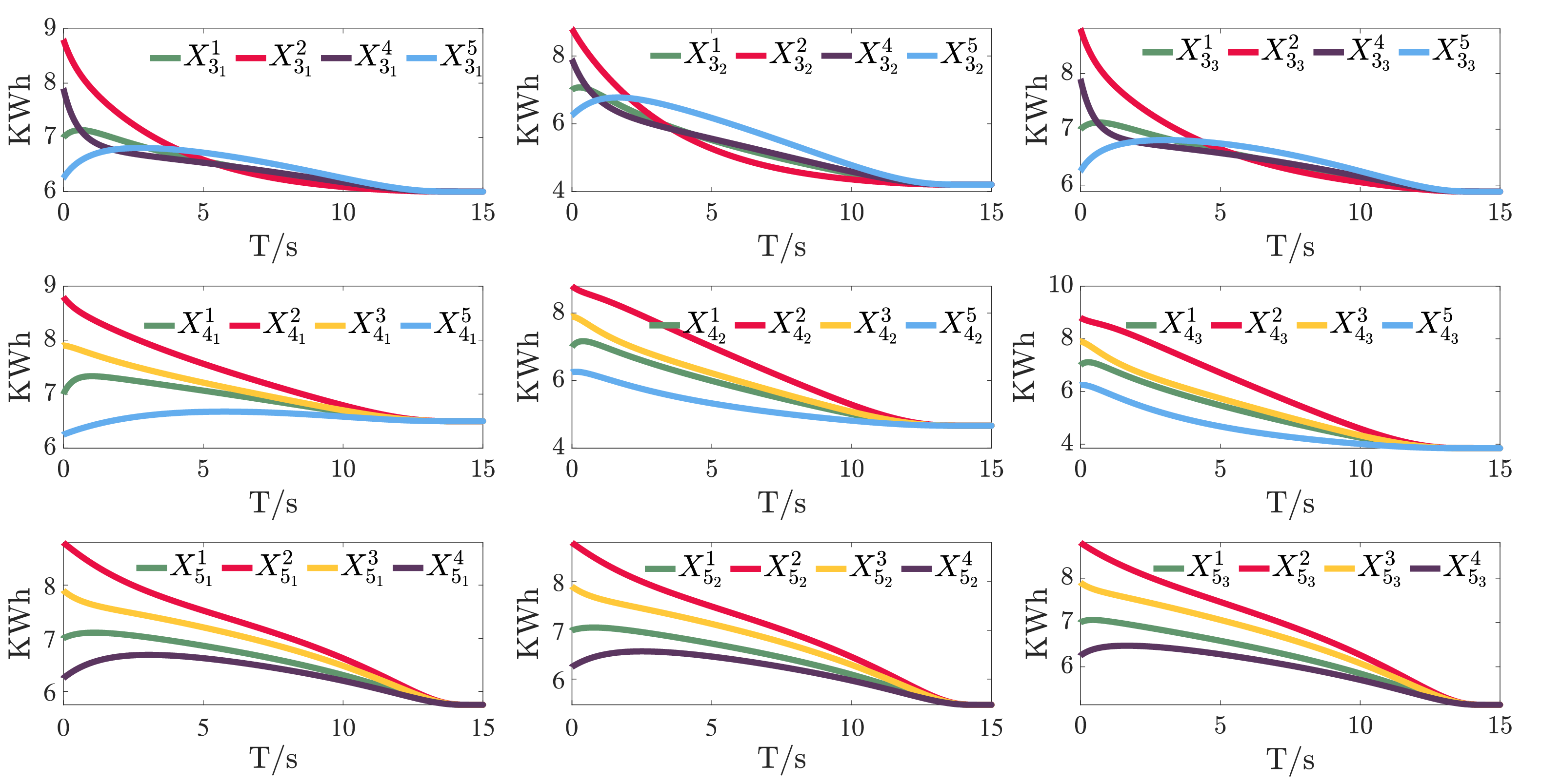}
\caption{Diagrams of State Variables Varying with Time}
\label{state}
\end{figure}

As is shown in Fig. \ref{state}, the electricity consumption for three grades appliances reach a consensus at the prescribed-time $T_{p}$. Thus, the theoretical result is verified.
\section{Conclusion}\label{con}
This paper considers adaptive fully distributed penalty-based NE seeking problem for noncooperative game with private inequality constraints and prescribed-time performance. A novel NE seeking algorithm is proposed based on adaptive penalty technique with a time-varying control gain and a cosine-based function. Compared with previous literature, this paper proposes an algorithm with faster convergence rate, more distribution features and lower computation burden.
\bibliographystyle{IEEEtran}
\bibliography{reference.bib}

\begin{thebibliography}{10}
\providecommand{\url}[1]{#1}
\csname url@samestyle\endcsname
\providecommand{\newblock}{\relax}
\providecommand{\bibinfo}[2]{#2}
\providecommand{\BIBentrySTDinterwordspacing}{\spaceskip=0pt\relax}
\providecommand{\BIBentryALTinterwordstretchfactor}{4}
\providecommand{\BIBentryALTinterwordspacing}{\spaceskip=\fontdimen2\font plus
\BIBentryALTinterwordstretchfactor\fontdimen3\font minus
  \fontdimen4\font\relax}
\providecommand{\BIBforeignlanguage}[2]{{%
\expandafter\ifx\csname l@#1\endcsname\relax
\typeout{** WARNING: IEEEtran.bst: No hyphenation pattern has been}%
\typeout{** loaded for the language `#1'. Using the pattern for}%
\typeout{** the default language instead.}%
\else
\language=\csname l@#1\endcsname
\fi
#2}}
\providecommand{\BIBdecl}{\relax}
\BIBdecl

\bibitem{Nash1950}
J.~F. Nash \emph{et~al.}, ``Non-cooperative games,'' 1950.

\bibitem{Maskery2009}
M.~Maskery, V.~Krishnamurthy, and Q.~Zhao, ``Decentralized dynamic spectrum
  access for cognitive radios: Cooperative design of a non-cooperative game,''
  \emph{IEEE Transactions on Communications}, vol.~57, no.~2, pp. 459--469,
  2009.

\bibitem{Ershen2024}
W.~Ershen, L.~Fan, H.~Chen, G.~Jing, Z.~Lin, X.~Jian, and H.~Ning,
  ``{MADRL}-based {UAV} swarm non-cooperative game under incomplete
  information,'' \emph{Chinese Journal of Aeronautics}, 2024.

\bibitem{Sun2020}
C.~Sun and G.~Hu, ``Continuous-time penalty methods for {N}ash equilibrium
  seeking of a nonsmooth generalized noncooperative game,'' \emph{IEEE
  Transactions on Automatic Control}, vol.~66, no.~10, pp. 4895--4902, 2020.

\bibitem{Sun2021}
------, ``Distributed generalized {N}ash equilibrium seeking for monotone
  generalized noncooperative games by a regularized penalized dynamical
  system,'' \emph{IEEE Transactions on Cybernetics}, vol.~51, no.~11, pp.
  5532--5545, 2021.

\bibitem{Jia2021}
W.~Jia, N.~Liu, and S.~Qin, ``An adaptive continuous-time algorithm for
  nonsmooth convex resource allocation optimization,'' \emph{IEEE Transactions
  on Automatic Control}, vol.~67, no.~11, pp. 6038--6044, 2021.

\bibitem{Ye2017}
M.~Ye and G.~Hu, ``Distributed nash equilibrium seeking by a consensus based
  approach,'' \emph{IEEE Transactions on Automatic Control}, vol.~62, no.~9,
  pp. 4811--4818, 2017.

\bibitem{Ye2019}
M.~Ye, G.~Hu, F.~L. Lewis, and L.~Xie, ``A unified strategy for solution
  seeking in graphical {$ N $}-coalition noncooperative games,'' \emph{IEEE
  Transactions on Automatic Control}, vol.~64, no.~11, pp. 4645--4652, 2019.

\bibitem{Gadjov2018}
D.~Gadjov and L.~Pavel, ``A passivity-based approach to nash equilibrium
  seeking over networks,'' \emph{IEEE Transactions on Automatic Control},
  vol.~64, no.~3, pp. 1077--1092, 2018.

\bibitem{Bianchi2020}
M.~Bianchi and S.~Grammatico, ``Fully distributed {N}ash equilibrium seeking
  over time-varying communication networks with linear convergence rate,''
  \emph{IEEE Control Systems Letters}, vol.~5, no.~2, pp. 499--504, 2020.

\bibitem{Ye2021}
M.~Ye and G.~Hu, ``Adaptive approaches for fully distributed {N}ash equilibrium
  seeking in networked games,'' \emph{Automatica}, vol. 129, p. 109661, 2021.

\bibitem{Ma2022}
T.~Ma, Z.~Deng, and C.~Hu, ``A fully distributed {N}ash equilibrium seeking
  algorithm for n-coalition games of {E}uler--{L}agrange players,'' \emph{IEEE
  Transactions on Control of Network Systems}, vol.~10, no.~1, pp. 205--213,
  2022.

\bibitem{Tao2022}
Q.~Tao, Y.~Liu, C.~Xian, and Y.~Zhao, ``Prescribed-time distributed
  time-varying {N}ash equilibrium seeking for formation placement control,''
  \emph{IEEE Transactions on Circuits and Systems II: Express Briefs}, vol.~69,
  no.~11, pp. 4423--4427, 2022.

\bibitem{Feng2020}
Z.~Feng and G.~Hu, ``Prescribed-time fully distributed {N}ash equilibrium
  seeking in noncooperative games,'' \emph{arXiv preprint arXiv:2009.11649},
  2020.

\bibitem{Zhao2023}
Y.~Zhao, Q.~Tao, C.~Xian, Z.~Li, and Z.~Duan, ``Prescribed-time distributed
  {N}ash equilibrium seeking for noncooperation games,'' \emph{Automatica},
  vol. 151, p. 110933, 2023.

\bibitem{Sun2023}
C.~Sun and G.~Hu, ``{N}ash equilibrium seeking with prescribed performance,''
  \emph{Control Theory and Technology}, vol.~21, no.~3, pp. 437--447, 2023.

\bibitem{Song2023}
Y.~Song, H.~Ye, and F.~L. Lewis, ``Prescribed-time control and its latest
  developments,'' \emph{IEEE Transactions on Systems, Man, and Cybernetics:
  Systems}, vol.~53, no.~7, pp. 4102--4116, 2023.

\bibitem{Gong2021}
X.~Gong, Y.~Cui, J.~Shen, J.~Xiong, and T.~Huang, ``Distributed optimization in
  prescribed-time: {T}heory and experiment,'' \emph{IEEE Transactions on
  Network Science and Engineering}, vol.~9, no.~2, pp. 564--576, 2021.

\bibitem{Zhu2023}
Q.~Zhu, \emph{Nonlinear systems}.\hskip 1em plus 0.5em minus 0.4em\relax
  MDPI-Multidisciplinary Digital Publishing Institute, 2023.

\bibitem{Cui2021}
L.~Cui and N.~Jin, ``Prescribed-time {ESO}-based prescribed-time control and
  its application to partial igc design,'' \emph{Nonlinear Dynamics}, vol. 106,
  no.~1, pp. 491--508, 2021.

\end{thebibliography}
\newpage

\end{document}